\documentclass[11pt]{article}
\usepackage{amsfonts,latexsym}
\input xy
\xyoption{all}
\newcommand{\cleqn}{\setcounter{equation}{0}}
\newcommand{\clth}{\setcounter{theorem}{0}}
\newcommand {\sectionnew}[1]{\section{#1}\cleqn\clth}
\newcommand{\beq}{\begin{equation}}
\newcommand{\eeq}{\end{equation}}
\newcommand{\beqa}{\begin{eqnarray}}
\newcommand{\eeqa}{\end{eqnarray}}
\newcommand{\beaa}{\begin{eqnarray*}}
\newcommand{\ben}{\begin{eqnarray*}}
\newcommand{\eaa}{\end{eqnarray*}}
\newcommand{\een}{\end{eqnarray*}}

\newcommand{\text}{\textrm}
\newcommand \nc {\newcommand}
\nc \proof {\noindent {\em{Proof.\/ }}}
\nc \qed {$\Box$\hfill}
\newtheorem{theorem}{Theorem}[section]
\newtheorem{lemma}[theorem]{Lemma}
\newtheorem{proposition}[theorem]{Proposition}
\newtheorem{corollary}[theorem]{Corollary}
\newtheorem{definition}[theorem]{Definition}
\newtheorem{example}[theorem]{Example}
\newtheorem{remark}[theorem]{Remark}
\newtheorem{conjecture}[theorem]{Conjecture}
\newtheorem{question}[theorem]{Question}
\nc \bth[1] {\begin{theorem}\label{t#1} }
\nc \ble[1] {\begin{lemma}\label{l#1} }
\nc \bpr[1] {\begin{proposition}\label{p#1} }
\nc \bco[1] {\begin{corollary}\label{c#1} }
\nc \bde[1] {\begin{definition}\label{d#1}\rm }
\nc \bex[1] {\begin{example}\label{e#1}\rm }
\nc \bre[1] {\begin{remark}\label{r#1}\rm }
\nc \bcon[1] {\begin{conjecture}\label{con#1}\rm }
\nc \bque[1] {\begin{question}\label{que#1}\rm }
\nc {\eth} { \end{theorem} }
\nc {\ele} { \end{lemma} }
\nc {\epr}{ \end{proposition} }
\nc {\eco} { \end{corollary} }
\nc {\ede} {\end{definition} }
\nc {\eex} { \end{example} }
\nc {\ere} {\end{remark} }
\nc {\econ} { \end{conjecture} }
\nc {\eque} {\end{question} }
\nc \eqref[1] {{\rm{(\ref{#1})}}}
\nc \thref[1]{Theorem \ref{t#1}}
\nc \leref[1]{Lemma \ref{l#1}}
\nc \prref[1]{Proposition
\ref{p#1}} \nc \coref[1]{Corollary \ref{c#1}}
\nc \deref[1]{Definition \ref{d#1}}
\nc \exref[1]{Example \ref{e#1}}
\nc \reref[1]{Remark \ref{r#1}}
\nc \conref[1]{Conjecture\ref{con#1}}

\def \R {{\mathcal R}}

\def\a{\alpha}
\def\b{\beta}

\def\Ga{\Gamma}

\def\l{\lambda}

\def\d{\partial}

\def \Z {{\mathbb Z}}
\def \Q {{\mathbb Q}}
\def \R {{\mathbb R}}

\def \tr{ {\mathrm{Tr}}}

\def \Im { {\mathrm{Im}} }

\renewcommand \ker { {\mathrm{Ker}} }

\nc \Wr {Wr} \nc \GRN { \Gr^{(N)} }
\nc \GRA[1] { \Gr_A^{(#1)} }   
\nc \GRAN { \GRA{N} } \nc \GrA[1] { \Gr_A(#1) }\nc \GrAa {
\GrA{\alpha} }
\nc \GRB[1] { \Gr_B^{(#1)} }   
\nc \GRBN { \GRB{N} } \nc \GrB[1] { \Gr_B(#1) } \nc \GrBb {
\GrB{\beta} }
\nc \GRMB[1] { \Gr_{MB}^{(#1)} }   
\nc \GRMBN { \GRMB{N} } \nc \GrMB[1] { \Gr_{MB}(#1) } \nc \GrMBb {
\GrMB{\beta} }

\begin{document}
\title{{\LARGE\bf{Cohomology of $GL_4(\Z)$ with Non-trivial
Coefficients}}}

\author{
I. ~Horozov
\thanks{E-mail: ihorozov@brandeis.edu}
\\ \hfill\\ \normalsize \textit{Department of Mathematics,}\\
\normalsize \textit{ Brandeis University, 415 South St.,}\\
\normalsize \textit {MS 050, Waltham, MA 02454 }  \\
}
\date{}

\maketitle

\begin{abstract}
In this paper we compute the cohomology groups of $GL_4(\Z)$ with
coefficients in symmetric powers of the standard representation
twisted by the determinant. This problem arises in Goncharov's
approach to the study of motivic multiple zeta values of depth 4.
The techniques that we use include Kostant's formula for
cohomology groups of nilpotent Lie subalgebras of a reductive Lie
algebra, Borel-Serre compactification, a result of Harder on
Eisenstein cohomology. Finally, we need to show that the ghost
class, which is present in the cohomology of the boundary of the
Borel-Serre compactification, disappears in the Eisenstein
cohomology of $GL_4(\Z)$. For this we use a computationally
effective version for the homological Euler characteristic of
$GL_4(\Z)$ with non-trivial coefficients.
\end{abstract}
\tableofcontents

\section{Introduction}
\subsection{Main result and applications}

The main goal of this paper is to present a computation of
cohomology groups
$$H^i(GL_4(\Z), S^{n-4} V_4 \otimes det),$$
where $S^{n-4} V_4$ is the $(n-4)$-th
symmetric power of the standard representation $V_4$
and $det$ is the determinant representation.

The above cohomology groups describe certain spaces
of motivic multiple zeta values.
This relation was revealed by Goncharov who suggested to me the problem
of computing the cohomology groups of $GL_4(\Z)$.

Recall the definition
of multiple zeta values
$$\zeta(k_1,\dots,k_m)=
\sum_{0<n_1<\dots <n_m}\frac{1}{n_1^{k_1}\dots n_m^{k_m}},$$
where $k_1+\dots +k_m$ is called {\it{weight}} and $m$ is
called {\it{depth}}.

Goncharov has described the cases of {\it{depth}}=$2$ \cite{G2}
and of {\it{depth}}=$3$ \cite{G3}. He relates the space of motivic
multiple zeta values of {\it{depth}}=$2$ and {\it{weight}}=$n$ to
the cohomology groups of $GL_2(\Z)$ with coefficients in the
$(n-2)$-symmetric power of the standard representation $V_2$,
namely, to
$$H^i(GL_2(\Z),S^{n-2}V_2).$$
He calls this a misterious relation between the multiple zeta
values of {\it{depth}}=$m$ and the ''modular variety''
$$GL_m(\Z)\backslash GL_m(\R)/SO_m(\R)\times \R^{\times}_{>0}.$$
In the paper \cite{G3}, he relates the spaces of
motivic multiple zeta values
of {\it{depth}}=$3$ and {\it{weight}}=$n$ to the cohomology of
$GL_3(\Z)$ with
coefficients in the $(n-3)$-symmetric power of
the standard representation $V_3$, namely,
$$H^i(GL_3(\Z),S^{n-3}V_3).$$

Goncharov has also related
the case of multiple zeta values of
{\it{depth}}=$4$ and {\it{weight}}=$n$ to the computation of
the cohomology of $GL_4(\Z)$ with coefficients in the
$(n-4)$-symmetric power
of the standard representation $V_4$ twisted by the determinant
(private communications). That is,
in order to compute the spaces of
motivic multiple zeta values of {\it{depth}}=$4$
and {\it{weight}}=$n$ one has to know
$$H^i(GL_4(\Z),S^{n-4}V_4 \otimes det).$$

The main result of this paper is the following. \bth{1.1} The
dimensions of the cohomology groups of $GL_4(\Z)$ with
coefficients the symmetric powers of the standard representation
twisted by the determinant are given by

$H^i(GL_4(\Z),S^{n-4}V_4 \otimes det)=
\left\{\begin{tabular}{ll}
$\Q \oplus H^1_{cusp}(GL_2(\Z),S^{n-2}V_2 \otimes det)$ &for $i=3,$\\
$0$                                    &for $i\neq 3.$
\end{tabular}\right.$
\\
\\
More explicitly,

$dim(H^3(GL_4(\Z),S^{12n-4+k}V_4 \otimes det))
=\left\{\begin{tabular}{ll}
$n+1$ & for $k=0,4,6,8,10,$\\
$n$   & for $k=2,$\\
$0$   & for $k$ odd.
\end{tabular}\right.$
\eth

\subsection{Computational methods and notation}
All representations that we consider are finite dimensional
representations of $GL(\Q)$ defined over $\Q$. However, we shall
consider them as representations of the arithmetic subgroups via
inclusion. We assume that the reader is familiar with group
cohomology. For a good introduction to this subject and to various
Euler characteristics of group see \cite{Br}.

We are going to describe briefly various types of cohomology groups of
arithmetic groups, namely, boundary cohomology, cohomology at the
infinity, Eisenstein cohomology, interior cohomology and cusp
cohomology. All of them are based on a compactification of certain
space, called Borel-Serre compactification. The reader who is not
familiar with these constructions should not be discouraged. We
have tried to present a piece of ''Calculus'' for cohomology of
arithmetic groups. That is, we give the definitions intuitively
rather than strictly, and describe the computational tools which
we are going to use. The constructions and the proofs of the basic
tools could be found in the cited literature. What we do in the
main part of this paper is to present the desired computation
based on these tools.

We start with the Borel-Serre compactification
\cite{BoSe}. Let $\Ga$ be a subgroup of $GL_m(\Q)$ which is
commensurable
to $GL_m(\Z)$. That is, the intersection $\Ga \cap GL_m(\Z)$ is
of finite index both in $\Ga$ and in $GL_m(\Z)$.
Let
$$X=GL_m(\R)/SO_m(\R)\times\R^{\times}_{>0}.$$
Then $X$ is a contractable topological space on which $\Ga$ acts on the left.
And let
$$Y_{\Ga}=\Ga\backslash X.$$
Then the Borel-Serre compactification of $Y_\Ga$, denoted by
$\overline{Y}_\Ga$,
is a compact space, containing $Y_\Ga$. Moreover, it is of the
same homotopy
type as $Y_\Ga$. If $V$ is a representation of $\Ga$ and
$V^\sim$ is the corresponding sheaf then
$$H^i_{top}(\overline{Y}_\Ga,V^{\sim})=H^i_{group}(\Ga,V).$$

The space $\overline{Y}_\Ga$ can be split into
strata, where each stratum corresponds
to a parabolic subgroup $P$ of $GL_{m/\Q}$
and the maximal stratum is $Y_\Ga$.
Also the closure of a stratum corresponding
to a parabolic subgroup $P$ consists of all strata corresponding
 to parabolic
subgroups $Q$ so that $Q\subset P$.
Let $Y_{\Ga,P}$ be the stratum corresponding to a
parabolic subgroup $P$. Let
$$P(\Z)=P(\Q)\cap \Ga.$$
Then the topological
cohomology of  $\overline{Y}_P$ coincides with the group cohomology of
$P(\Z)$. More precisely,
$$H^i_{top}(\overline{Y}_P,j^*_P V^{\sim})=H^i_{group}(P(\Z),V),$$
where $V$ is a representation over the rational numbers and
$V^{\sim}$ the corresponding sheaf on $\overline{Y}_\Ga$ and
$j_P^*V^\sim$ is its restriction on $\overline{Y}_{\Ga,P}$.

The boundary of the Borel-Serre compactification is
$$\d\overline{Y}_\Ga=\overline{Y}_\Ga-Y_\Ga=\cup_P \overline{Y}_{\Ga,P}.$$
The inclusion $$j:\d\overline{Y}_\Ga\subset \overline{Y}_\Ga$$ induces
$$j^{\#}: H^i_{top}(\overline{Y}_\Ga,V^{\sim})\rightarrow
H^i_{top}(\d\overline{Y}_\Ga,j^*V^{\sim}).$$
We call the range of the last map $j^\#$ cohomology of the boundary.
We use the notation
$$H^i_\d(\Ga,V):=H^i_{top}(\d\overline{Y}_\Ga,j^*V^{\sim}).$$
We warn the reader that it is not a standard notation.

The image of the map $j^\#$ is called cohomology at the infinity
of $\Ga$. We use the notation
$$H^i_{inf}(\Ga,V):=\Im(j^\#).$$
And the kernel of the map $j^\# $ is called interior cohomology of $\Ga$.
We use the notation
$$H^i_{!}(\Ga,V):=\ker(j^\#).$$

For the representations that we will consider we have that the
cohomology at infinity coincides with the Eisenstein cohomology.
This is used for describing certain maps between cohomology
groups. Also the interior cohomology coincides with the cusp
cohomology. In the representations which we will concider we are
going to use that fact in order to show that the interior cohomology
vanishes.

In our problem we have
$$H^i_{cusp}(GL_4(\Z),S^{n-4}V_4\otimes det)=0,$$
where $V_m$ is the standard $m$-dimensional representation of
$GL_m(\Q)$. And $S^{n}$ is the $n$-th symmetric power. The last
equality holds for $n>4$ because the representation
$$S^{n-4}V_4\otimes det$$
is not
self-dual. For $n=4$ it is true because $$H^i_{cusp}(SL_4(\Z),\Q)=0.$$
Thus, we need to compute only the Eisenstein cohomology.

The highest weight representation will be denoted by $L[a_1, ...
,a_m]$, where the weight $[a_1, ... ,a_m]$ sends $diag[H_1, ...
,H_m]$ to $a_1(H_1)+ ... +a_m(H_m)$. Sometimes we shall denote the
weight simply by $\l$. At a later stage there will be a number of
cohomologies to consider. In order to make the answer more
observable, sometimes we abbreviate. For example:
$$H^i(L[a_1,\dots,a_d]):=H^i(GL_d(\Z), L[a_1,\dots,a_d]).$$
For further abbreviation we set
$$\begin{array}{ll}
&(a_1,a_2|a_3|a_4):=H^1(L[a_1,a_2]) \otimes H^0(L[a_3]) \otimes H^0(L[a_4])\\
&(a_1|a_2,a_3|a_4):=H^0(L[a_1]) \otimes H^1(L[a_2,a_3])  \otimes H^0(L[a_4])\\
&(a_1|a_2|a_3,a_4):= H^0(L[a_1]) \otimes H^0(L[a_2]) \otimes H^1(L[a_3,a_4])\\
&(a_1,a_2|a_3,a_4):=H^1(L[a_1,a_2]) \otimes H^1(L[a_3,a_4])\\
&(a_1|a_2|a_3|a_4):=H^0(L[a_1]) \otimes H^0(L[a_2]) \otimes H^0(L[a_3])
\otimes H^0(L[a_4])
\end{array}$$
We also will use the abbreviation
$$(\overline{a_1,a_2}|a_3|a_4):=H^1_{cusp}(L[a_1,a_2]) \otimes H^0(L[a_3])
\otimes H^0(L[a_4]).$$ We consider the parabolic subgroups of
$GL_4$ that contain a fixed Borel subgroup. We shall consider the
standard representation of
 $GL_4$ with the choice of the Borel subgroup $B$ being the upper
triangular matrices. Then the parabolic subgroups can be listed in
the following way: $P_{ij}$ is the smallest parabolic subgroup
containing a non-zero $a_{ji}$-entry. And $P_{12,34}$ is the
smallest parabolic subgroup containing $a_{21} \neq 0$ and $a_{43}
\neq 0$. More precisely: All parabolic subgroups contain $B$ which
is upper triangular. Also, $P_{12}$ has a quotient $GL_2 \times
GL_1 \times GL_1,$ $P_{23}$ has a quotient $GL_1 \times GL_2
\times GL_1,$ $P_{34}$ has a quotient $GL_1 \times GL_1 \times
GL_2,$ $P_{13}$ has a quotient $GL_3 \times GL_1,$ $P_{24}$ has a
quotient $GL_1 \times GL_3,$ and $P_{12,34}$ has a quotient $GL_2
\times GL_2.$

We are going to use the Kostant's theorem \cite{K} in order to
obtain information about the parabolic subgroups. To do that we
need to examine carefully the action of the Weyl group, $W$ on the
root system of $gl_n$. Also we need the Weyl group, $W_P$
associated to the algebra $P$. In order to use the Kostant
theorem, we need to examine the action of the Weyl group $W$ on
the root system of $gl_n$ up to permutation of the root system of
$P$. That is, we need to consider representatives of the quotient
$W_P \backslash W$. We state Kostant's theorem \cite{K}.
\bth{1.1}  Let $V$ be a representation of highest
weight $\lambda$. Let $N_P$ be nilpotent radical of a parabolic group $P$,
and let $\rho$ be half of the sum of the positive roots. Then
$$H^i(N_P, V)= \oplus _{\omega} L_{{\omega}(\lambda + \rho) - \rho},$$
where the sum is taken over the representatives of the quotient
$W_P \backslash W$ with minimal length such that their length is
exactly $i$. In the above notation $L_\lambda$ means
representation of $N_P$ with highest weight $\lambda.$ \eth

Let $[a,b,c,d]$ denote an element of the root lattice (inside
$h^*$) whose value on the diagonal entry
$[H_{11},H_{22},H_{33},H_{44}]$ in $h$ is $a H_{11}+ b H_{22}+ c
H_{33}+ d H_{44}.$ The Weyl group acts on the weight lattice by
permuting the entries of $[a,b,c,d].$ It is well known that the
Weyl group is generated by reflections perpendicular to the
primitive roots. We can choose positivity so that the primitive
roots correspond to the permutation $(12)$, $(23)$ and $(34)$,
(having $sl_4$ in mind; $(12)$ sends $[a,b,c,d]$ to  $[b,a,c,d]$.)
Then the length of an element of the Weyl group is precisely the
(minimal) number of successive transpositions, or equivalently,
the (minimal) number of reflections w.r.t. the primitive roots. In
this setting the right quotient  $W_P \backslash W$ can be
interpreted as shuffles in the following way: Take for example the
parabolic subalgebra $P_{23}$. Its Levi quotient $M_P =
M_{P_{23}}$ is $gl_1 \times gl_2 \times gl_1.$ Thus, $W_P$ is
generated by $(23).$ Among the representatives of the quotient
$W_P \backslash W$ we can consider the ones that preserve the
order of the subset $\{23\}$ inside $\{1234\}$. Thus, we can
consider all shuffles of $\{1|23|4\}.$ Similarly, if we take the
parabolic subalgebra $P_{12,34}$, we need to consider the shuffles
of $\{12|34\}$ so that the order of $\{12\}$ and the order of $\{34\}$ is preserved. 
And for the subalgebras $P_{13}$ we consider the
shuffles of the set $\{123|4\}$, which means permutations of
$\{1234\}$ such that the order $\{123\}$ is preserved.

In order to apply Kostant's theorem, we need to examine the length of
each element $\omega$
in the Weyl group $W,$ and also the resulting weight
$\omega (\lambda + \rho) - \rho$, where
$\lambda=[a,b,c,d]$ is the weight of $V$ and
$\rho$ is half of the sum of the positive roots.

After we obtain the cohomology of the parabolic groups we have to
consider a spectral sequence involving these cohomologies in order
to obtain the cohomology of the boundary of the Borel-Serre
compactification. Then we use homological Euler characteristics in order to compute
the cohomology groups of $GL_m(\Z)$ for $m=2,3,4$.

{\bf Acknowledgments:} I would like to thank Professor Goncharov
for giving me this problem and for computational techniques that I
learned from him. I would like to thank Professor Harder for
teaching me important computational techniques.

This work was initiated at Max-Planck Institute f\"ur Mathematik.
I am very grateful for the stimulating atmosphere, created there,
as well as for the financial support during my stay.

\section {Homological Euler characteristics of $GL_m(\Z)$}
We call {\it{homological Euler characteristic of a group}} $\Gamma$
the alternating sum of the dimension of the cohomology of the group.
We denote it by
$\chi_h(\Gamma,V),$
where $V$ is a finite dimensional representation of $\Gamma$. More precisely,
$$\chi_h(\Gamma,V)=\sum_i (-1)^i dim H^i(\Gamma,V).$$
In this section we compute the homological Euler characteristics
of $GL_m(\Z)$ for $m=2,3,4$ with representations which later will
occur in the Kostant's formula applied to $GL_4(\Z)$ with
coefficients in the representation $(n-4)$-th symmetric power of
the standard representation twisted by the determinant which is
$L[n-3,1,1,1]$

The material in this section is in the spirit of the papers \cite{Ho2}
and \cite{Ho1}. Most of the formulas and notations are taken from there.
The only exception is the computation of $\chi_h(GL_3(\Z),L[n-3,1,0])$, done here in details.

We start with $GL_2(\Z)$.
\bth{3.1}
   Let $S^nV_2$ be the $n$-th symmetric power of the
standard representation of $GL_2$. Then\\

$\chi_h(GL_2(\Z), S^{12n+k}V_2)= \left\{\begin{tabular}{cl}
$-n+1$ & $k=0$ \\
$-n$    & $k=2,4,6,8$ \\
$-n-1$  & $k=10$ \\
$0$     & $k=odd,$ \\
\end{tabular}\right.
     $\\
            and\\

$\chi_h(GL_2(\Z), S^{12n+k}V_2\otimes \det)=
\left\{\begin{tabular}{cl}
$-n$   & $k=0$ \\
$-n-1$ & $k=2,4,6,8$ \\
$-n-2$ & $k=10$ \\
$0$    & $k=odd.$ \\
   \end{tabular}\right.
    $
        \\
     \eth

For $GL_m(\Z)$ $m=3 \mbox{ and }4$ we need to consider the representations\\

\begin{tabular}{l}
$L[n-3,1,0]=\ker(S^{n-3}V_3\otimes V_3 \rightarrow S^{n-2}V_3),$\\
\\
$L[n-2,1,1]=S^{n-3}V_3\otimes det,$\\
\\
$L[n-2,2,2]=S^{n-4}V_3,$ \\
\\
$L[n-3,1,1,1]=S^{n-4}V_4\otimes det.$\\
\end{tabular}
\bth{3.2}
The homological Euler characteristics of $GL_3(\Z)$ and $GL_4(\Z)$
with coefficients in the above representation are given by\\

\begin{tabular}{l}
(a) $\chi_h(GL_3(\Z),L[n-3,1,0])
 =\chi_h(GL_2(\Z),S^{n-4}V_2)-\chi_h(GL_2(\Z),S^{n-2}V_2)$,\\
\\
(b) $\chi_h(GL_3(\Z),L[n-2,1,1])=\chi_h(GL_2(\Z),S^{n-2}V_2\otimes det),$\\
\\
(c) $\chi_h(GL_3(\Z),L[n-2,2,2])=\chi_h(GL_2(\Z),S^{n-4}V_2),$ \\
\\
(d) $\chi_h(GL_4(\Z),L[n-3,1,1,1])=\chi_h(GL_2(\Z),S^{n-2}V_2\otimes det).$\\

\end{tabular}
\eth The technique that we are going to use involves a substantial
simplification of the trace formula which works when
$\Gamma=GL_m(\Z)$ or a group co-mensurable to $GL_m(\Z)$. The
simplification of the trace formula for $GL_m(\Z)$ was developed
in \cite{Ho1,Ho2}. Besides the simplification we are going to use
some computation which were done in the above two papers.

Now we present the simplification of the trace formula in the case of
$GL_m(\Z)$. An arithmetic group $\Gamma$
has also an orbifold Euler characteristic.
We denote it by $\chi(\Gamma)$, without subscript.
It is in fact an Euler characteristic of a certain orbifold.
There is
a more algebraic description. If an arithmetic group $\Gamma$ has no torsion
then the orbifold Euler characteristic coincides with the homological
Euler characteristic with coefficients in the trivial representation.
$$\chi(\Gamma)=\chi_h(\Gamma,\Q).$$
If $\Gamma$ has torsion choose a torsion free finite index subgroup $\Gamma_0$.
Then
$$\chi(\Gamma)=\frac{\chi(\Gamma_0)}{[\Gamma:\Gamma_0]}.$$

Let $C(A)$ denote the centralizer of the element $A$ inside $\Gamma$.
Then the classical trace formula is
$$\chi_h(\Gamma,V)=\sum_{A} \chi(C(A))\tr(A|V),$$
where the sum is taken over all torsion elements considered up to
conjugation. And $C(A)$ denotes the centralizer of the element $A$
inside $\Gamma$. We remark that in this formula the identity
element is also considered as a torsion element.

For the simplification of the trace formula we need the following
definition. Let $A$ be an element in $GL_m(\Z)$. Consider it as an
$m\times m$ matrix. Let $f$ be its characteristic polynomial. Let
$$f=f_1^{a_1}\dots f_l^{a_l}$$
be the factorization of $f$ into irreducible over $\Q$
polynomials. Denote by
$$R(g,h)=\prod_{i,j}(\a_i-\b_j)$$ the resultant of the polynomials
$$g=\prod_i(x-\a_i) \mbox{ and } h=\prod_j(x-\b_i).$$
Denote by
$$R(A)=\prod_{i<j}R(f_i^{a_i},f_j^{a_j})$$

\bth{2.10}
   Let $V$ be a finite dimensional representation
of $GL_m(\Q)$. Then the homological Euler characteristic of
$GL_m(\Z)$ with coefficients in $V$ is given by
 $$
   \chi_h(GL_m(\Z),V)=
\sum_{A} |R(A)|\chi(C(A))\tr(A|V),
    $$
  where the sum is taken over torsion matrices $A$ consisting
of square blocks $A_{11},\dots A_{ll}$ on the block-diagonal and zero
blocks off the diagonal. Also the matrices $A_{ii}$ are non-conjugate to
each other. And they are chosen from the set
$\{+1,+I_2,-1,-I_2,T_3,T_4,T_6\},$
where
$$T_3=
\left[
\begin{tabular}{rr}
$0$ & $1$\\
$-1$ & $-1$
\end{tabular}
\right],
\mbox{ }
T_4=
\left[
\begin{tabular}{rr}
$0$ & $1$\\
$-1$ & $0$
\end{tabular}
\right],
\mbox{ }
T_6=
\left[
\begin{tabular}{rr}
$0$ & $-1$\\
$1$ & $1$
\end{tabular}
\right].
$$
The blocks on the diagonal are chosen up to permutation. And the
characteristic polynomial $f_i$ of $A_{ii}$ is a power of an
irreducible polynomial, and $f_i$ and $f_j$ are relatively prime.
    \eth
{\bf{Remark:}} There is one more simplification that we can make.
In the formula in theorem 3.3 for the homological Euler characteristic
one can do the summation in the following way. If the $-I_m$ acts
 on $V$ nontrivially then all the cohomologies of $GL_m(\Z)$ with coefficients
in $V$ vanish and the homological Euler characteristic vanishes.
If $-I_m$ acts trivially on $V$ then $\tr(-A|V)=\tr(A|V)$. Also, $C(-A)=C(A)$
and $|R(-A)|=|R(A)|$. If $-A$ is not conjugate to $A$ then in the
sum of theorem 3.3 we can compute the invariants for $A$. And for $-A$
they are the same. Note that $-A$ is conjgate to $A$ if and only if
one can obtain $-A$ by permuting the blocks on the diagonal of $A$.

\proof (of theorem 3.2) Parts (b), (c) and (d) are computed \cite{H1}.
We are going to prove part (a). We are going to use the following notation.
Given a matrix $A$ whose blocks no the diagonal are $A_{11},\dots, A_{ll}$
and whose blocks off the diagonal are zero, we write it as
$$A=[A_{11},\dots, A_{ll}].$$ Using this notation and the notation of
theorem 2.3 we quote lemma 4.2 of the paper\cite{Ho1}
\ble{4.2}For the centralizers and the resultants of the torsion elements in
$GL_3(\Z)$ we have\\

\begin{tabular}{l}
(c) $|R([I_2,-1])|\chi(C([I_2,-1]))=-\frac{1}{12}$,\\
\\
(e) $|R([T_3,1])|\chi(C([T_3,1]))=\frac{1}{4}$,\\
\\
(f) $|R([T_6,1])|\chi(C([T_6,1]))=\frac{1}{12}$,\\
\\
(i) $|R([T_4,1])|\chi(C([T_4,1]))=\frac{1}{4}$.\\
\\
\end{tabular}
\ele
Also, we are going to use lemma 5.3 from the same paper \cite{H1}.
\ble{5.3}
    The traces of the torsion
elements in $GL_3(\Z)$ acting on the symmetric power of the standard
representation are given by:\\
\\
(c) \mbox{ } $\tr([I_2,-1]|S^{2n+k}V_3)= \left\{\begin{tabular}{cc}
$n+1$ & $k=0$ \\ $n+1$ & $k=1$, \\
\end{tabular}\right.
$ \\
(e) \mbox{ } $\tr([T_3,1]|S^{3n+k}V_3)=
\left\{\begin{tabular}{cc} $1$ & $k=0$ \\ $0$ & $k=1$ \\ $0$ &
$k=2$, \\
\end{tabular}\right.
$\\
(f) \mbox{ } $\tr([T_6,1]|S^{6n+k}V_3)
=
\left\{\begin{tabular}{rc} $1$ & $k=0$ \\ $2$ & $k=1$ \\ $2$ &
$k=2$ \\ $1$ & $k=3$ \\ $0$ & $k=4$ \\ $0$ & $k=5$, \\
\end{tabular}\right.
$\\
\\
(i)  \mbox{ } $\tr([T_4,1]|S^{4n+k}V_3)= \left\{\begin{tabular}{cc}
$1$  & $k=0$ \\ $1$  & $k=1$ \\ $0$  & $k=2$ \\ $0$  & $k=3$. \\
\end{tabular}\right.
$\\
\ele

In order to compute $\tr(A|L[w-3,1,0])$ for torsion elements $A$,
we are going to use
$$L[w-3,1,0]=\ker(S^{w-3}V_3\otimes V_3 \rightarrow S^{w-2}V_3).$$
Also, we are going to use that
$$\tr(A|V\otimes W)=\tr(A|V)\tr(A|W).$$
Using the above two equalities together with lemma 3.5, we obtain
the following.
\ble{5.3}
    The traces of the torsion
elements in $GL_3(\Z)$ acting on the $L[w-3,1,0]$ are given by:\\
\\
(c) \mbox{ } $\tr([I_2,-1]|L[2n-1+k,1,0])= \left\{\begin{tabular}{cc}
$-1$ & $k=0$ \\ $0$ & $k=1$, \\
\end{tabular}\right.
$ \\
(e) \mbox{ } $\tr([T_3,1]|L[3n-1+k,1,0])=
\left\{\begin{tabular}{cc} $-1$ & $k=0$ \\ $0$ & $k=1$ \\ $0$ &
$k=2$, \\
\end{tabular}\right.
$\\
(f) \mbox{ } $\tr([T_6,1]|L[6n-1+k,1,0])
=
\left\{\begin{tabular}{rc} $-1$ & $k=0$ \\ $0$ & $k=1$ \\ $2$ &
$k=2$ \\ $3$ & $k=3$ \\ $2$ & $k=4$ \\ $0$ & $k=5$, \\
\end{tabular}\right.
$\\
\\
(i)  \mbox{ } $\tr([T_4,1]|L[4n-1+k,1,0])= \left\{\begin{tabular}{cc}
$-1$  & $k=0$ \\ $0$  & $k=1$ \\ $1$  & $k=2$ \\ $0$  & $k=3$. \\
\end{tabular}\right.
$\\
\ele

For each of the torsion elements $A$ in $GL_3(\Z)$ we have that
$A$ and $-A$ are not conjugated. When we use theorem 2.3 we can
count only four of the torsion elements listed in lemmas 2.4, 2.5
and 2.6 and multiply by two in order to consider the contribution
of the negative of these torsion elements. Thus, using theorem
2.3, lemma 2.4 and
lemma 2.6 we obtain\\
\\
\begin{tabular}{l}
$\chi_h(GL_3(\Z),L[12n-1,1,0])=2(\frac{1}{12}-\frac{1}{4}-\frac{1}{12}
-\frac{1}{4})=-1,$\\
\\
$\chi_h(GL_3(\Z),L[12n+1,1,0])=2(\frac{1}{12}+0+\frac{2}{12}
+\frac{1}{4})=1,$\\
\\
$\chi_h(GL_3(\Z),L[12n+3,1,0])=2(\frac{1}{12}+0+\frac{2}{12}
-\frac{1}{4})=0,$\\
\\
$\chi_h(GL_3(\Z),L[12n+5,1,0])=2(\frac{1}{12}-\frac{1}{4}-\frac{1}{12}
+\frac{1}{4})=0,$\\
\\
$\chi_h(GL_3(\Z),L[12n+7,1,0])=2(\frac{1}{12}+0+\frac{2}{12}
-\frac{1}{4})=0,$\\
\\
$\chi_h(GL_3(\Z),L[12n+1,1,0])=2(\frac{1}{12}+0+\frac{2}{12}
+\frac{1}{4})=1,$\\
\end{tabular}\\

Consider the statement of theorem 2.2 part (a). The above computation
of homological Euler characteristics
gives the left hand side of part(a).
The right hand side can be computed directly
from theorem 2.1. They do coincide. Thus, part (a) of theorem 2.2 is
 proven.
\section {Cohomology of $GL_2(\Z)$.}
This section is to show how the computational method works for
$GL_2(\Z)$. All the results are known, but we need them for the
later sections. We are going to compute Eisenstein cohomology and cusp cohomology of
$GL_2(\Z)$ with coefficients in some representations. 

First we are going to compute the cohomology of the boundary using
Kostant's theorem.
Let $L[a,b]$ be the irreducible representation with
highest weight $[a,b]$. The group $GL_2$ has one parabolic subgroup up to conjugation - 
the Borel subgroup $B$.
It has a nilpotent radical $N$ and a Levi quotient $GL_1 \times GL_1$.
The Weyl group
has two elements. Also, the half of the `sum' of the positive roots
is $\rho=[1/2,-1/2]$ Consider the following table:


$$
\begin{array}{llllllll}
&\omega \epsilon W &length &\omega (\l +\rho )-\rho \\
&12                &0      &[a,b]\\
&21                &1      &[b-1,a+1]\\
\end{array}
$$

From Kostant's theorem we obtain that
$$H^n(N, L[a,b])= \left\{
\begin{array}{lll}
&L[a,b]     &n=0,\\
&L[b-1,a+1] &n=1.
\end{array}
\right.
$$
The integral points of the
Levi quotient of $B$ are $GL_1(\Z) \times GL_1(\Z)$. Using
the Hochschild-Serre spectral sequence we compute $H^n(B,L[a,b])$.
If both $a$ and $b$ are even then
$H^0(B,L[a,b])=H^0(GL_1(\Z), L[a]) \otimes H^0(GL_1(\Z), L[b])=\Q$,
and the rest of the cohomology groups are trivial.
If both $a$ and $b$ are odd then
$H^1(B,L[a,b])=H^0(GL_1(\Z), L[b-1]) \otimes H^0(GL_1(\Z), L[a+1])=\Q$.
If $a+b$ is odd then $H^n(B,L[a,b])=0$ for all $n$.

There are several cases. If $a+b$ is odd then $-I$ acts non-trivially
on $L[a,b]$. So the cohomology of $GL_2(\Z)$ vanishes. If $a=b=2k$ then
$L[a,b]$ is the trivial representation of $GL_2(\Z)$. So
$$H^i(GL_2(\Z),L[2k,2k])=H_{Eis}^i(GL_2(\Z),L[2k,2k])=
\left\{\begin{tabular}{ll}
$\Q$ & $i=0,$\\
$0$  & $i=1,$
\end{tabular}
\right.
$$
and
$$H^i_{cusp}(GL_2(\Z),L[2k,2k])=0.$$

If $a=b=2k+1$ then
$$H^i(GL_2(\Z),L[2k+1,2k+1])=0.$$
So the Eisenstein and the cusp cohomology also vanish.

The interesting cases are when both $a$ and $b$ are even or when
both $a$ and $b$ are odd. For those cases we do not give a
complete proof, but rather an interpretation of the cohomologies.
If follows from considering modular forms for $GL_2(\Z)$ of weight $2(a-b)$ or equivalently, holomorphic modular forms for $SL_2(\Z)$. The Eisenstein cohomology is generated by the Eisenstein series and the dimension of the cusp cohomology $H^1_{cusp}(GL_2(\Z),L[a,b])$ is equal to the dimension of the cusp forms of weight $2(a-b)$.
In any of these cases we have
$$H^0(GL_2(\Z),L[a,b])=0.$$ Also,
if $a$ and $b$ are both odd, we have that the map
$$H^1(GL_2(\Z),L[a,b])\rightarrow H^1(B,L[a,b])=\Q$$
is surjective. Then
$$H^1_{Eis}(GL_2(\Z),L[a,b])=\Q,$$
and
$$dimH^1_{cusp}(GL_2(\Z),L[a,b])=-1+dimH^1(GL_2(\Z),L[a,b]).$$
If the weights $a$ and $b$ are both even, then the Eisenstein cohomology
coincides with the whole group cohomology.

Here is one interpretation of the cohomology of $GL_2(\Z)$ in
cases when both $a$ and $b$ are either
even or odd.
We are not going to use the following interpretation, only the above formulas,
but it is nice to keep it in mind.

Let $a$ and $b$ be both odd.
Then
$$\begin{tabular}{ll}
&$H^1(SL_2(\Z),L[a,b])=$\\
&$H^1_{cusp}(GL_2(\Z),L_{[a+1,b+1]})\oplus H^1_{cusp}(GL_2(\Z),L_{[a,b]})
\oplus H_{Eis}^1(GL_2(\Z),L_{[a,b]}).$
\end{tabular}
$$
The first direct summand corresponds to holomorphic
cuspidal forms of weight $a-b-2$. The second summand correspond to
anti-holomorphic cusp forms of weight $a-b-2$. And the last summand
 corresponds to the Eisenstein series of weight $a-b-2$ (when bigger than 2).

Keeping in mind the above decompositions one can compute the dimensions of
the cohomlogy groups (or dimensions of cusp forms) using theorem 2.1.
Note that in theorem 2.1 the homological Euler characteristic is
equal to minus the dimension of the first cohomology group, since the higher
cohomology groups vanish as well as the zeroth.
\sectionnew{Cohomology of $GL_3(\Z)$}
In this section we compute cohomology groups of $GL_3(\Z)$ with coefficients in
certain representations which are needed for our main problem.
They arise as representations of the Levi quotients of two of the
maximal parabolic subgroups of $GL_4$, namely, $P_{13}$ and $P_{24}$.
They lead to computation of cohomology groups of $GL_3(\Z)$
 with coefficients
in any of the representations
$L[0,0,0]=\Q$, $L[w-3,1,0]$,
$L[w-2,2,2]$ and $L[w-2,1,1]$.
\bth{5.1} The cohomology of $GL_3(\Z)$ with coefficients
in the above representations are given by
\\
\\
\begin{tabular}{l}
(a) $H^i(GL_3(\Z), \Q)=\left\{
\begin{array}{lll}
&(0|0|0)  &i=0,\\
&0        &i \neq 0.
\end{array}
\right.$
\\
\\
(b) $H^i(GL_3(\Z),L[n-3,1,0])=\left\{
\begin{array}{lll}
&(\overline{n-3,-1}|2) &i=2\\
&(-2|\overline{n-2,2}) &i=3\\
&0                     &i\neq 2,3
\end{array}
\right.$
\\
\\
(c) $H^i(GL_3(\Z), L[n-2,2,2])=\left\{
\begin{array}{lll}
&(0|\overline{n-1,3}) &i=3\\
&0                    &i \neq3
\end{array}
\right.$
\\
\\
(d) $H^i(GL_3(\Z), L[n-2,1,1])=\left\{
\begin{array}{lll}
&(0|\overline{n-1,1}) &i=2\\
&0                    &i \neq 2.
\end{array}
\right.
$
\end{tabular}
\eth

Before proving the above theorem, we examine
the cohomology of $GL_3(\Z)$ with coefficients in $L_{[a,b,c]}$

The algebraic group $GL_3$ has three parabolic subgroups: $B$,
$P_{12}$ and $P_{23}$. In order to find their cohomology groups, we
need the explicit action of the Weyl group; more precisely we need
the various $\omega(\l+\rho)-\rho$ that enter in Kostant's
theorem. Note that half of the sum of the positive roots is
$\rho=[1,0,-1]$.



$$
\begin{array}{llllllll}
&\omega \epsilon W  &\mbox{length of }\omega& &\omega(\l) &\omega(\l+\rho)-\rho \\
&123                &\:\:0&            &[a,b,c]    &[a,b,c]\\
&132                &\:\:1&            &[a,c,b]    &[a,c-1,b+1]\\
&213                &\:\:1&            &[b,a,c]    &[b-1,a+1,c]\\
&231                &\:\:2&            &[b,c,a]    &[b-1,c-1,a+2]\\
&312                &\:\:2&            &[c,a,b]    &[c-2,a+1,b+1]\\
&321                &\:\:3&            &[c,b,a]    &[c-2,b,a+2]
\end{array}
$$

Using the Kostant's theorem we find the cohomology groups
of the nilpotent radicals of the parabolic groups.

$$
H^q(H,L[a,b,c])=\left\{
\begin{array}{lll}
&L[a,b,c]                             &q=0  \\
&L[a,c-1,b+1]   \oplus L[b-1,a+1,c]   &q=1  \\
&L[b-1,c-1,a+2] \oplus L[c-2,a+1,b+1] &q=2  \\
&L[c-2,b,a+2]                         &q=3
\end{array}
\right.
$$

$$
H^q(N_{12},L[a,b,c])=\left\{
\begin{array}{lll}
&L[a,b,c]        &q=0\\
&L[a,c-1,b+1]    &q=1\\
&L[b-1,c-1,a+2]  &q=2
\end{array}
\right.
$$

$$
H^q(N_{23},L[a,b,c])=\left\{
\begin{array}{lll}
&L[a,b,c]        &q=0\\
&L[b-1,a+1,c]    &q=1\\
&L[c-2,a+1,b+1]  &q=2
\end{array}
\right.
$$

In order to pass to cohomologies of the parabolic groups,
we use the Hochschild-Serre spectral sequence relating
the nil radical and the Levi quotient of a parabolic subgroup
to the parabolic subgroup itself; namely the short exact
sequence $N \rightarrow P \rightarrow S$. We recall the notation
$H^n(L[a_1,...,a_k])=H^n(GL_k \Z, L[a_1,...,a_k])$ and
$(a|b|c)= H^0(L[a]) \otimes H^0(L[b]) \otimes H^0(L[c])$.

$$
H^i(B,L[a,b,c])=\left\{
\begin{array}{lll}
&(a|b|c)                            &i=0\\
&(a|c-1|b+1)   \oplus (b-1|a+1|c)   &i=1\\
&(b-1|c-1|a+2) \oplus (c-2|a+1|b+1) &i=2\\
&(c-2|b|a+2)                        &i=3
\end{array}
\right.
$$

$$
E_2^{p,q}(P_{12},L[a,b,c])=\left\{
\begin{array}{lll}
&H^p(L[a,b])    \otimes H^0(L[c])   &q=0\\
&H^p(L[a,c-1])  \otimes H^0(L[b+1]  &q=1\\
&H^p(L[b-1,c-1])\otimes H^0(L[a+2]  &q=2
\end{array}
\right.
$$

$$
E_2^{p,q}(P_{23},L[a,b,c])=\left\{
\begin{array}{lll}
&H^0(L[a])   \otimes H^p(L[b,c])     &q=0\\
&H^0(L[b-1]  \otimes H^p(L[a+1,c])   &q=1\\
&H^0(L[c-2]) \otimes H^p(L[a+1,b+1]) &q=2
\end{array}
\right.
$$

It is true that the above two spectral sequences
stabilize at the $E_2$-level.
However, in any particular case the formulas will be
much simpler, and one can use them to compute
the boundary cohomology.

Let $B$, $P_{12}$ $P_{23}$ be the parabolic subgroups of $GL_3 \Z$.

{\bf $H^i(GL_3(\Z),\Q)$}

For part (a) we have
$$
H^i(B, \Q)=\left\{
\begin{array}{lll}
&(0|0|0)  &i=0\\
&(-2|0|2) &i=3\\
&0        &n \neq 0,3
\end{array}
\right.
$$

$$H^0(P_{12}, \Q)=H^0(GL_2(\Z) ,\Q) \otimes H^0(GL_1(\Z), \Q)$$

$$H^0(P_{23}, \Q)=H^0(GL_1(\Z) ,\Q) \otimes H^0(GL_2(\Z), \Q)$$

From Mayer-Vietoris we obtain that the boundary cohomology of $GL_3(\Z)$ is

$$H^i_{\partial}(GL_3(\Z),\Q)=\left\{
\begin{array}{lll}
&(0|0|0)  &i=0\\
&(-2|0|2) &i=4\\
&0        &i \neq 0,4
\end{array}
\right.
$$

The homological Euler characteristic of $GL_3(\Z)$ with
trivial coefficients is $1$ (see theorem 2.2 part (c) and theorem 2.1).
That is,
$$\chi_h(GL_3(\Z),\Q)=1.$$

Then the forth cohomology of the boundary component disappears in
the Eisenstein cohomology. Therefore,

$$H_{Eis}^i(GL_3(\Z),\Q)=\left\{
\begin{array}{lll}
&(0|0|0)  &i=0\\
&0        &i \neq 0.
\end{array}
\right.
$$
Also, the cusp cohomology of $GL_3(\Z)$ with trivial coefficients is zero.
Therefore the Eisenstein cohomology coincides with the whole group cohomology.

We proceed to part(b).

{\bf $H^i(B,L[n-3,1,0])$}

Using the computations in the beginning of this section, we obtain

$$H^i(B,L[n-3,1,0])=\left\{
\begin{array}{lll}
&0          &i=0\\
&(0|n-2|0)  &i=1\\
&(-2|n-2|2) &i=2\\
&0          &i=3
\end{array}
\right.
$$

$$H^i(P_{12},L[n-3,1,0])=\left\{
\begin{array}{lll}
&(n-3,1|0)  &i=1\\
&(n-3,-1|2) &i=2\\
&0          &i\neq 1,2
\end{array}
\right.
$$

$$H^i(P_{23},L[w-3,1,0])=\left\{
\begin{array}{lll}
&(0|n-2,0)  &i=2\\
&(-2|n-2,2) &i=3\\
&0          &i\neq2,3\\
\end{array}
\right.
$$

Using Mayer-Vietoris, for the cohomology of the boundary
of the Borel-Serre compactification, we obtain

$$H^i_\partial(GL_3(\Z),L[n-3,1,0])=\left\{
\begin{array}{lll}
&(\overline{n-3,1}|0)                             &i=1\\
&(\overline{n-3,-1}|2)\oplus (0|\overline{n-2,0}) &i=2\\
&(-2|\overline{n-2,2})                            &i=3\\
&0                                                &i\neq 1,2,3
\end{array}
\right.
$$

The representation $L[n-3,1,0]$ is not self dual. So the cohomology of
$GL_3(\Z)$ with coefficients in $L[n-3,1,0]$ coincides with the
Eisenstein cohomology, which is a subspace of the cohomology
of the boundary. The first cohomology of $GL_3(\Z)$ with
coefficients in any representation vanishes. For the homological
Euler characteristic of $GL_3(\Z)$ with coefficients in $L[n-3,1,0]$
(theorem 2.2 part (a)) we have
$$\chi_h(GL_3(\Z),L[n-3,1,0])=
\chi_h(GL_2(\Z),S^{n-4}V_2)-\chi_h(GL_2(\Z),S^{n-2}V_2).$$
We obtain that the dimension of the second cohomology is half of the
dimension of the second cohomology of the boundary of the
Borel-Serre compactification. That is,

$$dim H^2_{Eis}(GL_3(\Z),L[n-3,1,0])=
\frac{1}{2} dim H^2_\partial(GL_3(\Z),L[n-3,1,0]).$$
Also,
$$dim H^3_{Eis}(GL_3(\Z),L[n-3,1,0])=
 dim H^3_\partial(GL_3(\Z),L[n-3,1,0]).$$
The second cohomology of the boundary is a direct sum of two spaces with
the same dimesions. In order to find out which of the subspaces or
which linear combination of the spaces enters in the Eisenstein cohomology,
we have to consider the central characters of the two parabolic subgroups \cite{Ha}. 
For the parabolic subgroup $P_{12}$ we take the central torus
$$\left[
\begin{tabular}{ccc}
$t$ &     &         \\
    & $t$ &         \\
    &     & $t^{-2}$
\end{tabular}
\right].$$
The highest weight induces a character on it, namely $[n-3,-1,2]$,
whose evaluation on the above
element is
$$n-3-1-2\times 2=n-8.$$
For the parabolic subgroup $P_{23}$ we take the central torus
$$\left[
\begin{tabular}{ccc}
$t^2$ &          &         \\
      & $t^{-1}$ &         \\
      &          & $t^{-1}$
\end{tabular}
\right].$$
The highest weight induces a character on it, namely $[0,n-2,0]$,
whose evaluation on the above
element is
$$0-(n-2)=-n+2.$$
Their sum is -6. The space which enters in the Eisenstein cohomology
has higher weight. Thus we need to solve
$$n-8>-n+2.$$
Thus for $n>5$ we have
$$H^i(GL_3(\Z),L[n-3,1,0])=\left\{
\begin{array}{lll}
(\overline{n-3,-1}|2) &i=2\\
(-2|\overline{n-2,2}) &i=3\\
0                    &i\neq 2,3
\end{array}
\right.
$$
The value of $n$ is always even and greater or equal to $4$. The other option
for $n$ is $n=4$. Then
$$H^i(GL_3(\Z),L[1,1,0])=\left\{
\begin{array}{lll}
(0|\overline{4-2,0}) &n=2\\
(-2|\overline{4-2,2}) &n=3\\
0                    &n\neq 2,3
\end{array}
\right.
$$
That is,
$$H^i(GL_3(\Z),L[1,1,0])=0.$$

{\bf $H^*(GL_3, L[n-2,2,2])$ when $n$ is even.}

Using the computation in the beginning of section 3 we obtain:

$$H^i(B,L[n-2,2,2])=\left\{
\begin{array}{lll}
&(n-2|2|2) &i=0\\
&0         &i=1\\
&0         &i=2\\
&(0|2|n)   &i=3
\end{array}
\right.
$$

$$H^i(P_{12},L[n-2,2,2])=\left\{
\begin{array}{lll}
&(n-2,2|2) &i=1\\
&0         &i \neq1
\end{array}
\right.
$$

$$H^i(P_{23},L[n-2,2,2])=\left\{
\begin{array}{lll}
&(n-2|2|2]) &i=0\\
&(0|n-1,3)  &i=3\\
&0          &i \neq 0,3
\end{array}
\right.
$$

Using Mayer-Vietoris we obtain
$$H^i_\d(GL_3(\Z), L[n-2,2,2])=\left\{
\begin{array}{lll}
&(\overline{n-2,2}|2) &i=1\\
&(0|\overline{n-1,3}) &i=3\\
&0                    &i \neq 1,3
\end{array}
\right.
$$
The first cohomology of $GL_3(\Z)$ vanishes, Therefore,
$$H^i(GL_3, L[n-2,2,2])=\left\{
\begin{array}{lll}
&(0|\overline{n-1,3}) &i=3\\
&0                    &i \neq3
\end{array}
\right.
$$

{\bf $H^*(GL_3, L[n-2,1,1])$ when $n$ is even.}

Using the computation in the beginning of section 3 we obtain:

$$H^i(B,L[n-2,1,1])=\left\{
\begin{array}{lll}
&0         &i=0\\
&(n-2|0|2) &i=1\\
&(0|0|n)   &i=2\\
&0         &i=3
\end{array}
\right.
$$

$$H^i(P_{12},L[n-2,1,1])=\left\{
\begin{array}{lll}
&(n-2,0|2) \oplus (0|0|n) &i=2\\
&0                        &i \neq 2
\end{array}
\right.
$$

$$H^i(P_{23},L[n-2,1,1])=\left\{
\begin{array}{lll}
&(0|n-1,1)  &i=2\\
&0          &i \neq 2
\end{array}
\right.
$$

Using Mayer-Vietoris we obtain

$$H^i_{\partial}(GL_3(\Z), L[n-2,1,1])=\left\{
\begin{array}{lll}
&(0|0|n)\oplus (n-2,0|2) \oplus (0|\overline{n-1,1}) &i=2\\
&0                                                   &i \neq 2,
\end{array}
\right.
$$
 From the homological
Euler characteristic of $GL_3(\Z)$ with coefficients in $L[n-2,1,1]$
we obtain that
$$dim H^2_{Eis}(GL_3(\Z),L[n-2,1,1])=
\frac{1}{2}(-1+ dim H^2_\partial(GL_3(\Z),L[n-3,1,0])).$$
For the parabolic subgroup $P_{12}$ we take the central torus
$$\left[
\begin{tabular}{ccc}
$t$ &     &         \\
    & $t$ &         \\
    &     & $t^{-2}$
\end{tabular}
\right].$$
The highest weight induces a character on it, namely $[n-2,0,2]$,
whose evaluation on the above
element is
$$n-2+0-2\times 2=n-6.$$
For the parabolic subgroup $P_{23}$ we take the central torus
$$\left[
\begin{tabular}{ccc}
$t^2$ &          &         \\
      & $t^{-1}$ &         \\
      &          & $t^{-1}$
\end{tabular}
\right].$$
The highest weight induces a character on it, namely $[0,n-1,1]$,
whose evaluation on the above
element is
$$0-(n-1)-1=-n.$$
Their sum is -6. The space which enters in the Eisenstein cohomology
has higher weight. Thus we need to solve
$$n-6>-n.$$
Thus for $n>3$, which is always the case, we have
$$H^i(GL_3(\Z),L[n-2,1,1])=\left\{
\begin{tabular}{ll}
$(n-2,0|2)\oplus (0|0|n)$ & $i=2$\\
$0$                       & $i\neq 2$
\end{tabular}
\right.
$$
\section  {Cohomologies of the parabolic subgroups of $GL_4$.}
This section consists of computation of cohomology of the
parabolic subgroups of $GL_4(\Z)$ with coefficients in the
representation $S^{n-4} V_4 \otimes det$. We use Kostant's theorem
in order to compute these cohomology groups. In the process we reduce
the question to computation of the cohomology groups of the Levi
quotients which have factors $GL_1(\Z)$, $GL_2(\Z)$ or/and
$GL_3(\Z)$. For the last three groups we use the computation from
the sections on cohomolgy of $GL_2(\Z)$ and of $GL_3(\Z)$.

Recall the notation of the parabolic subgroups: We choose the
Borel subgroup $B$ to be the group of upper triangular matrices.
Let $N$ be its unipotent radical of $B$. Let $P_{ij}$ be the
smallest parabolic subgroup containing $B$ and containing a non
zero $a_{ji}$-entry. Similarly, $P_{12,34}$ is the smallest
(parabolic) subgroup containing $B$ and containing non zero
$a_{21}$- and $a_{43}$-entries. The unipotent radicals of $P_{ij}$
will be denoted by $N_{ij}$; and the Levi quotient by
$S_{ij}=P_{ij}/N_{ij}$.
\bpr{4.1}(cohomologies of the parabolic subgroups)
Let
\\
$V=S^{n-4} V_4 \otimes det$.
Then

$$H^i(B, V)= \left\{
\begin{array}{llllll}
&(0|n-2|0|2)                    &i=2 \\
&(0|0|0|n)  \oplus (-2|w-2|2|2) &i=3 \\
&(-2|0|2|n)                     &i=6\\
&0                              &i \neq 2,3,6.
\end{array}
\right.$$

$$H^i(P_{12}, V)= \left\{
\begin{array}{llllll}
&(n-3, 1|0|2)                  &i=2, \\
&(0|0|0|w) \oplus (n-3,-1|2|2) &i=3, \\
&0                             &i \neq 2,3.
\end{array}
\right.$$

$$H^i(P_{23},V)=
\left\{\begin{array}{llllll}
&(0|n-2,0|2) \oplus (0|0,0|n) &i=3\\
&(-2|n-2,2|2)                 &i=4, \\
&0                            &i \neq 3,4.
\end{array}
\right.$$

$$H^n(P_{34},V)= \left\{
\begin{array}{llllll}
&(0|0|n-1,1)\oplus (-2|n-2|2|2) &i=3, \\
&(-2|0|n-1,3)                   &i=6 \\
&0                              &i \neq 3,6.
\end{array}
\right.$$

$$H^i(P_{13}, V)= \left\{
\begin{array}{llllll}
&(0|0|0|n)\oplus (\overline{n-3,-1}|2|2) &i=3, \\
&(-2|\overline{n-2,2}|2)                 &i=4, \\
&0                                       &i\neq 3,4.
\end{array}
\right.$$

$$H^i(P_{12,34}, V)= \left\{
\begin{array}{llllll}
&(n-3,-1|2|2) \oplus (0|0|n-1,1) &i=3 \\
&0                               &i \neq 3.
\end{array}
\right.$$

$$H^i(P_{24}, V)= \left\{
\begin{array}{llllll}
&(0|0|0|n)\oplus (0|n-2,0|2)  &i=3, \\
&(-2|0|\overline{n-1,3})      &i=6, \\
&0                            &i\neq 3,6..
\end{array}
\right.$$

\epr

The main tool in the proof will be Kostant's theorem
and Hochschild-Serre spectral sequence.
In terms of weights representation $S^{n-4} \otimes det$
is $L[n-3,1,1,1]$. We shall denote the representation $L[n-3,1,1,1]$ simply
by $V$. We identify the Weyl group of $GL_4$ with
the permutation group of four elements.
We also need the length of the permutation
which we denote by $l$.

$$
\begin{array}{llllllllllll}
\\
&w    \: &l \: &w(\lambda+\rho)-\rho \\
&1234 \: &0 \: &[w-3,1,1,1]\\
&1243 \: &1 \: &[w-3,1,0,2]\\
&1324 \: &1 \: &[w-3,0,2,1]\\
&1342 \: &2 \: &[w-3,0,0,3]\\
&1423 \: &2 \: &[w-3,-1,2,2]\\
&1432 \: &3 \: &[w-3,-1,1,3]\\
\\
&2134 \: &1 \: &[0,w-2,1,1]\\
&2143 \: &2 \: &[0,w-2,0,2]\\
&2314 \: &2 \: &[0,0,w-1,1]\\
&2341 \: &3 \: &[0,0,0,w]\\
&2413 \: &3 \: &[0,-1,w-1,2]\\
&2431 \: &4 \: &[0,-1,1,w]\\
\\
&3124 \: &2 \: &[-1,w-2,2,1]\\
&3142 \: &3 \: &[-1,w-2,0,3]\\
&3214 \: &3 \: &[-1,1,w-1,1]\\
&3241 \: &4 \: &[-1,1,0,w]\\
&3412 \: &4 \: &[-1,-1,w-1,3]\\
&3421 \: &5 \: &[-1,-1,2,w]\\
\\
&4123 \: &3 \: &[-2,w-2,2,2]\\
&4132 \: &4 \: &[-2,w-2,1,3]\\
&4213 \: &4 \: &[-2,1,w-1,2]\\
&4231 \: &5 \: &[-2,1,1,w]\\
&4312 \: &5 \: &[-2,0,w-1,3]\\
&4321 \: &6 \: &[-2,0,2,w]\\
\end{array}
$$

Now we can consider particular parabolic subgroup $P$.
In order to apply Kostan's theorem
we need to find good representatives
$W_P \backslash W$; more precisely representatives
of minimal length. This can be done by choosing the elements
of the permutation group
that preserve the ordered subsets corresponding to $P_{ij}$.
For example, when we consider $P_{23}$ the minimal representatives
of $W_{P_{23}}\backslash W$ are the permutations $w$ such that $w(2)<w(3)$.
When we consider $P_{12,34}$ we need the permutations $w$ such that
$w(1)<w(2)$ and $w(3)<w(4)$. And for the group $P_{13}$ the needed permutations
are the ones such that $w(1)<w(2)<w(3)$.
Thus, using Kostant's theorem we obtain:

It is easier to describe the cohomology
$$H^n(N,V),$$
than to write it down. One can think of it in the following way.
Consider the last column of the above table. If it is with
weight $[a,b,c,d]$ and with length $l$ then
$H^l(N,V)$ contains the representation
$L[a,b,c,d]$. Also,
all components of the cohomology are obtained in this way.

$$H^i(N_{12}, V)= \left\{
\begin{array}{llllll}
&L_{[n-3, 1, 1, 1]},                                                  &i=0,\\
&L_{[n-3, 1,0,+2]} \oplus L_{ [n-3,0,+2,1]},                          &i=1,\\
&L_{ [n-3,0,0,+3]} \oplus L_{ [n-3,-1,+2,+2]} \oplus L_{[0,0,n-1,1]}, &i=2,\\
&L_{[n-3,-1,1,+3]}  \oplus L_{[0,-1,n-1,2]} \oplus L_{ [0,0,0,n]},    &i=3,\\
&L_{[-1,-1,n-1,+3]} \oplus L_{[0,-1,1,n]},                            &i=4,\\
&L_{[-1,-1,2,n]},                                                     &i=5.
\end{array}
\right.
$$

$$H^i(N_{23}, V)= \left\{
\begin{array}{llllll}
&L_{[n-3, 1, 1, 1]},                                               &i=0, \\
&L_{[n-3, 1,0,+2]} \oplus L_{[0, n-2, 1, 1]},                      &i=1, \\
&L_{[n-3,0,0,3]} \oplus L_{ [0,n-2,0,+2]} \oplus L_{[-1,n-2,2,1]}, &i=2, \\
&L_{[-2,n-2,2,2]} \oplus L_{[-1,n-2,0,3]}\oplus L_{[0,0,0,n]},     &i=3, \\
&L_{[-2,n-2,1,3]} \oplus L_{[-1,1,0,n]},                           &i=4, \\
&L_{[-2,1,1,n]},                                                   &i=5.
\end{array}
\right.
$$

$$H^i(N_{34},V)= \left\{
\begin{array}{llllll}
&L_{[n-3, 1, 1, 1]},                                               &i=0, \\
&L_{[n-3,0,+2,1]}  \oplus L_{[0, n-2, 1, 1]},                      &i=1, \\
&L_{[n-3,-1,2,2]} \oplus L_{[-1,n-2,2,1]}\oplus L_{[0,0,n-1,1]},   &i=2, \\
&L_{[-2,n-2,2,2]} \oplus L_{[0,-1,n-1,2]}\oplus L_{ [-1,1,n-1,1]}, &i=3, \\
&L_{[-2,1,n-1,2]}\oplus L_{[-1,-1,n-1,+3]},                        &i=4, \\
&L_{[-2,0,n-1,3]},                                                 &i=5.
\end{array}
\right.
$$

$$H^i(N_{13}, V)= \left\{
\begin{array}{llllll}
&L_{[n-3, 1, 1, 1]} , &i=0, \\
&L_{[n-3, 1,0,+2]} ,  &i=1, \\
&L_{[n-3,0,0,3]} ,    &i=2, \\
&L_{[0,0,0,n]} ,      &i=3.
\end{array}
\right.
$$

$$H^i(N_{12,34},V)= \left\{
\begin{array}{llllll}
&L_{[n-3, 1, 1, 1]}                       &i=0 \\
&L_{[n-3,0,+2,1]}                         &i=1 \\
&L_{[n-3, -1,2,2]} \oplus L_{[0,0,n-1,1]} &i=2 \\
&L_{[0,-1,n-1,2]}                         &i=3 \\
&L_{[-1,-1,n-1,+3]}                       &i=4
\end{array}
\right.
$$

$$H^i(N_{24}, V)= \left\{
\begin{array}{llllll}
&L_{[n-3, 1, 1, 1]},  &i=0, \\
&L_{[0, n-2, 1, 1]},  &i=1, \\
&L_{[-1,n-2,2,1]},    &i=2, \\
&L_{[-2,n-2,2,2]},    &i=3.
\end{array}
\right.
$$

Now we apply the Hochschild-Serre spectral sequence to the exact sequences
$0 \rightarrow N_{ij} \rightarrow P_{ij} \rightarrow S_{ij} \rightarrow 0$.
Thus, the spectral sequence is
of the form $$H^p(S,H^q(N,V)) => H^{p+q}(P,V).$$
In the computation we are going
to use the Kunneth formula
$H^*(G_1 \times G_2, V_1 \otimes V_2) = H^*(G_1, V_1) \otimes H^*(G_2, V_2)$.
A substantial simplification comes from the facts that
$H^p(GL_m (\Z), det)=0$, for $m=1,2,3$.
It can be proven by the Hochschild-Serre spectral sequence
relating $GL_n$ to $SL_n$ and $G_m$.
One more observation about the computation of
the cohomology of the parabolic subgroups.
The cohomology groups of the nilpotent radical $H^*(N, V)$ are
representations of the Levi quotient.
For example, $H^0(N_{12}, S^{w-4}V_4 \otimes det)
= L_{[w-3, 1, 1, 1]} = L_{[w-3,1]} \otimes L_{[1]} \otimes L_{[1]}$,
since the Levi quotient is
$S_{12}=GL_2(\Z) \times GL_1(\Z) \times GL_1(\Z)$.

We are going to use some abbreviation in the computation that follows.
More precisely,
by $H^p(L[a,b])$ we mean  $H^p(GL_2 \Z, L_{[a,b]})$, similarly, by
$H^p(L[a])$ we mean  $H^p(GL_1 \Z, L_{[a]})$
and by $H^p(L[a,b,c])$ we mean  $H^p(GL_3 \Z, L_{[a,b,c]})$. Also,
we set $V = S^{w-4}V \otimes det=L[n-3,1,1,1]$.

\subsection{Cohomology of $B$}
The Levi quotient of a Borel subgroup is a Cartan subgroup.
Thus the representations
obtained from Kostan's theorem decompose into
tensor product of one dimensional representations.

$$
E_2 ^{p,q}=H^p(S, H^q(N,V)= \left\{
\begin{array}{llllll}
&H^p(S, L_{[0,n-2,0,2]})                               &q=2 \\
&H^p(S, L_{[0,0,0,n]}) \oplus H^p(S, L_{[-2,n-2,2,2]}) &q=3 \\
&H^p(S, L_{[-2,0,2,n]})                                &q=6\\
&0                                                     &q \neq 2,3,6.
\end{array}
\right.
$$

All other representations of $S$ do not contribute to the
cohomology of the Borel subgroup because at least one of the
entries of the weight is an odd number. The ones that are left
contain only even coefficients. Thus, they are trivial
representations of $GL_1 (\Z)$. Then the $E_2$-terms of the spectral
sequence can be simplified to:

$$
E_2 ^{p,q}= \left\{
\begin{array}{llllll}
&(0|n-2|0|2)                    &p=0, q=2 \\
&(0|0|0|n)  \oplus (-2|n-2|2|2) &p=0, q=3 \\
&(-2|0|2|n)                     &p=0, q=6\\
&0                              &otherwise
\end{array}
\right.
$$

The only non-zero entries of the above spectral sequence occur
only when $p=0$. Therefore the sequence degenerates at the
$E_2$-level, and cohomology of the Borel subgroup is

$$
H^i(B, V)= \left\{
\begin{array}{llllll}
&(0|n-2|0|2)                    &i=2 \\
&(0|0|0|n)  \oplus (-2|n-2|2|2) &i=3 \\
&(-2|0|2|n)                     &i=6\\
&0                              &i \neq 2,3,6
\end{array}
\right.
$$

\subsection{Cohomology of $P_{12}$}
We proceed similarly with the other parabolic subgroups. Recall,
the Levi quotient of $P_{12}$ is $S_{12}=GL_2(\Z) \times GL_1(\Z)
\times GL_1(\Z)$. Thus, the spectral sequence becomes:
$$E_2 ^{p,q}=H^p(S_{12}, H^q(N_{12}, V)$$
$$E_2 ^{p,q}=\left\{
\begin{array}{llllll}
&H^p(S_{12}, L_{[n-3, 1,0,+2]})   = H^p(S_{12}, L_{[n-3, 1]} \otimes L_{[0]}
\otimes L_{[2]}) &q=1, \\
&H^p(S_{12}, L_{ [n-3,-1,+2,+2]}) =H^p(S_{12}, L_{[n-3,-1]} \otimes L_{[2]}
\otimes L_{[2]})  &q=2, \\
&H^p(S_{12}, L_{ [0,0,0,n]})      = H^p(S_{12}, L_{[0,0]} \otimes L_{[0]}
\otimes L_{[n]})    &q=3, \\
&0                  &q \neq 1,2,3.
\end{array}
\right.
$$

The representations of the $GL_1(\Z)$ quotients that give a contribution
to the cohomology groups are the trivial representations.
Thus,

$$E_2 ^{p,q}= \left\{
\begin{array}{llllll}
&H^p(L[w-3, 1]) \otimes H^0(L[0]) \otimes H^0(L[2]) &q=1, \\
&H^p(L[w-3,-1])  \otimes H^0(L[2])^{\otimes 2}      &q=2, \\
&H^p(L[0,0]) \otimes H^0(L[0]) \otimes H^0(L[w])    &q=3, \\
&0                                                  &q \neq 1,2,3.
\end{array}
\right.
$$

We are going to use the fact that $H^p(GL_2(\Z), L)=0$ for $p>1$,
for any representation $L$ of $GL_2(\Q)$.
In particular the differential
$d_2: E_2 ^{p,q} \rightarrow E_2 ^{p+2,q-1}$ is zero, since
$E_2 ^{p,q}$ is non-zero only when $p=0$ or $p=1$.
Therefore, the spectral sequence degenerates at
the $E_2$-level, and the cohomology of $P_{12}$ is the following.

$$H^i(P_{12}, V)= \left\{
\begin{array}{llllll}
&(n-3, 1|0|2)                  &i=2, \\
&(0|0|0|n) \oplus (n-3,-1|2|2) &i=3, \\
&0                             &i \neq 2,3.
\end{array}
\right.
$$
\subsection{Cohomology of $P_{23}$}
Recall that the Levi quotient
$S_{23}$ of $P_{23}$ is $GL_1(\Z) \times GL_2(\Z) \times GL_1(\Z)$.
Then

$$E_2 ^{p,q}=H^p(S_{23}, H^q(N_{23},V)= \left\{
\begin{array}{llllll}
&H^p(S_{23}, L_{[0,n-2,0,2]})                        &q=2, \\
&H^p(S_{23}, L_{[-2,n-2,2,2]}  \oplus L_{[0,0,0,n]}) &q=3, \\
&0                                                   &q \neq 2,3.
\end{array}
\right.
$$

Using similar arguments, we obtain
$$H^i(P_{23},V)=
\left\{\begin{array}{llllll}
&(0|n-2,0|2) \oplus (0|0|0|n) &i=3\\
&(-2|n-2,2|2)                 &i=4, \\
&0                            &i \neq 3,4.
\end{array}
\right.
$$

\subsection{Cohomology of $P_{34}$}
Recall that the Levi quotient $S_{34}$ of $P_{34}$ is
$GL_1(\Z) \times GL_1(\Z) \times GL_2(\Z)$.
Then

$$E_2 ^{p,q}=H^p(S_{34}, H^q(N_{34},V))= \left\{
\begin{array}{llllll}
&H^p(S_{34},L_{[0,0,n-1,1]})  &q=2, \\
&H^p(S_{34},L_{[-2,n-2,2,2]}) &q=3, \\
&H^p(S_{34},L_{[-2,0,n-1,3]}) &q=5.
&0                            &q \neq 2,3,5
\end{array}
\right.
$$

Similarly, we obtain

$$E_2 ^{p,q} =\left\{
\begin{array}{llllll}
&H^0(L[0]) \otimes H^0(L[0]) \otimes H^p(L[n-1,1])  &q=2,\\
&H^0(L[-2]) \otimes H^0(L[n-2]) \otimes H^p(L[2,2]) &q=3,\\
&H^0(L[-2]) \otimes H^0(L[0]) \otimes H^p(L[n-1,3]) &q=5,\\
&0                                                  &q \neq 2,3,5
\end{array}
\right.
$$

And finally, the cohomology of $P_{34}$ is

$$H^i(P_{34},V)= \left\{
\begin{array}{llllll}
&(0|0|n-1,1)\oplus (-2|n-2|2|2) &i=3, \\
&(-2|0|n-1,3)                   &i=6 \\
&0                              &i \neq 3,6.
\end{array}
\right.
$$

\subsection{Cohomology of $P_{13}$}

Recall that the Levi quotient $S_{13}$ of $P_{13}$ is $GL_3 \Z \times GL_1 \Z$.

$$H^p(S_{13}, H^i(N_{13}, V))= \left\{
\begin{array}{llllll}
&0 &q=0, \\
&H^p(S_{13}, L_{[n-3, 1,0,2]})  &q=1, \\
&0 &q=2, \\
&H^p(S_{13}, L_{[0,0,0,n]})  &q=3.
\end{array}
\right.
$$

We can simplify it to

$$H^p(S_{13}, H^q(N_{13},V))= \left\{
\begin{array}{llllll}
&0                                  &q=0, \\
&H^p(L[n-3, 1,0])\otimes H^0(L[2])  &q=1, \\
&0                                  &q=2, \\
&H^p(L[0,0,0])\otimes H^0(L[n])     &q=3.
\end{array}
\right.
$$

From the section "Cohomology of $GL_3(\Z)$" we know that
for $n>5$ we have
$$H^p(GL_3(\Z),L[n-3,1,0])=\left\{
\begin{array}{lll}
(\overline{n-3,-1}|2) &p=2\\
(-2|\overline{n-2,2}) &p=3\\
0                     &p\neq 2,3
\end{array}
\right.
$$
And for $n=4$
$$H^p(GL_3(\Z),L[n-3,1,0])=0.$$
Also
$$H^p(GL_3(\Z),\Q)=\left\{
\begin{array}{lll}
&(0|0|0)  &p=0\\
&0        &p \neq 0.
\end{array}
\right.
$$
Therefore
$$H^p(S_{13}, H^q(N_{13},V))= \left\{
\begin{array}{llllll}
&(0|0|0|n)                &p=0\mbox{ and }q=3, \\
&(\overline{n-3,-1}|2|2)  &p=2\mbox{ and }q=1, \\
&(-2|\overline{n-2,2}|2)  &p=3\mbox{ and }q=1, \\
&0 &\mbox{for all other p and q}.
\end{array}
\right.
$$
The above spectral sequence degenerates at $E_2$ level. Therefore
$$H^i(P_{13},V))= \left\{
\begin{array}{llllll}
&(0|0|0|n) \oplus (\overline{n-3,-1}|2|2)  &i=3,\\
&(-2|\overline{n-2,2}|2)                   &i=4,\\
&0                                         &i\neq 3,4.
\end{array}
\right.
$$
\subsection{Cohomology of $P_{12,34}$}
For the last parabolic subgroup we can obtain a better answer
in terms of cohomology of $GL_2(\Z)$.
However, the $d_2$ differential might be non-trivial.
Recall that the Levi quotient $S_{12,34}$ of $P_{12,34}$ is
$GL_2(\Z) \times GL_2(\Z)$.
Then the spectral sequence is
$$E_2 ^{p,q}=H^p(S_{12,34},H^q(N_{12,34},V))$$
$$E_2 ^{p,q}=\left\{
\begin{array}{llllll}
&H^p(S_{12,34}, L_{[n-3, -1,2,2]})\oplus H^p(S_{12,34}, L_{[0,0,n-1,1]})&q=2\\
&0                                                                    &q\neq 2
\end{array}
\right.
$$
Therefore
$$E_2 ^{1,2}=
[H^1(L[n-3,-1])\otimes H^0(L[2,2])]
\oplus [H^0(L[0,0])\otimes H^1(L[n-1,1])].$$
And
$$E_2 ^{p,q}=0\mbox{ for } p\neq 1 \mbox{ or }q\neq 2$$
Finally,
$$H^i(P_{12,34},V)= \left\{
\begin{array}{llllll}
&(n-3,-1|2|2)\oplus (0|0|n-1,1) &i=3 \\
&0                              &i \neq 3
\end{array}
\right.
$$
\subsection{Cohomology of $P_{24}$}

Recall that the Levi quotient $S_{24}$ of $P_{24}$ is
$GL_1(\Z) \times GL_3(\Z)$.
We have the spectral sequence

$$H^p(S_{24}, H^q(N_{24},V))= \left\{
\begin{array}{llllll}
&0  &q=0, \\
&H^p(S_{24},L_{[0, n-2, 1, 1]}) &q=1, \\
&0   &q=2, \\
&H^p(S_{24},L_{[-2,n-2,2,2]}) &q=3.
\end{array}
\right.
$$

We can simplify it to

$$H^p(S_{24}, H^q(N_{24},V))= \left\{
\begin{array}{llllll}
&0  &q=0, \\
&H^0(L[0])\otimes H^p(L[n-2, 1, 1]) &q=1, \\
&0   &q=2, \\
&H^0(L[-2])\otimes H^p(L[n-2,2,2])  &q=3.
\end{array}
\right.
$$
From the section "Cohomology of $GL_3(\Z)$" we know that
$$H^i(GL_3(\Z),L[n-2,1,1])=\left\{
\begin{tabular}{ll}
$(n-2,0|2)\oplus (0|0|n)$ & $i=2$\\
$0$                       & $i\neq 2$
\end{tabular}
\right.
$$
And also
$$H^i(GL_3, L[n-2,2,2])=\left\{
\begin{array}{lll}
&(0|\overline{n-1,3}) &i=3\\
&0                    &i \neq3
\end{array}
\right.
$$
Therefore,
$$H^p(S_{24}, H^q(N_{24},V))= \left\{
\begin{tabular}{llllll}
$(0|0|0|n)\oplus (0|n-2,0|2)$  &$p=2,$ $q=1$, \\
$(-2|0|\overline{n-1,3})$      &$p=3,$ $q=3$, \\
$0$                            & for all other p and q.
\end{tabular}
\right.
$$
The spectral sequence degenerates. Therefore,
$$H^i(P_{24},V))= \left\{
\begin{tabular}{llllll}
$(0|0|0|n)\oplus (0|n-2,0|2)$  &$i=3,$ \\
$(-2|0|\overline{n-1,3})$      &$i=6,$ \\
$0$                            &$i\neq 3,6.$
\end{tabular}
\right.
$$
\sectionnew{Boundary cohomology of $GL_4(\Z)$}
In this section we compute the cohomology of the boundary of the
Borel-Serre compactification associated to $GL_4(\Z)$ with coefficients in
$$V=S^{n-4}\otimes det = L[n-3,1,1,1].$$ The Eisenstein cohomology,
which in our case is the whole group cohomology, injects
into the cohomology of the boundary.

We recall briefly several statements about Borel-Serre compactification
associated to $GL_m(\Z)$. Let
$$X=GL_m(\R)/SO_m(\R)\times\R^{\times}_{>0}.$$
And let
$$Y=GL_m(\Z)\backslash X.$$
Then the Borel-Serre compactification of $Y$, denoted by $\overline{Y}$,
is a compact space, containing $Y$, and of the same homotopy type.
The space $\overline{Y}$ is obtained by attaching cell $\sigma_P$ to $X$,
corresponding to each parabolic subgroup $P$. Denote by $Y_P$ the projection
of $\sigma_P$ to $\overline{Y}$. Let $\overline{Y}_P$ be the closure of
 $Y_P$. Then $\overline{Y}_Q \subset\overline{Y}_P$ when $Q\subset P$.
The boundary of $\overline{Y}$ is obtained by gluing together
the spaces $\overline{Y}_P$.
In the following computation we shall denote by $Y_{ij}$ the space
$Y_{P_{ij}}$. For these spaces we have
$$H^i_{top}(\overline{Y}_{ij},i^*F_V)=
H^i_{group}(P_{ij},V),$$
for a suitable sheaf $F_V$ on $\overline{Y}$, where $i$ is the inclusion
of $\overline{Y}_{ij}$ into $\overline{Y}$. For simplification we will
 not write the restriction functor $i^*$.

The cohomology of the boundary can be computed the spectral sequence
of the type 'Mayer-Vietoris'.

$$
\xymatrix{
   &H^q(\overline{Y}_{13},F_V))
    \ar[r] \ar[dr]
             &H^q(\overline{Y}_{12},F_V))  \ar[dr]\\
E_1^{*.q}: &H^q(\overline{Y}_{12,34},F_V)  \ar[ur]
       |!{[u];[r]}\hole \ar[dr] |!{[d];[r]}\hole
             & H^q(\overline{Y}_{23},F_V) \ar[r]
                    &  H^q(\overline{Y}_B,F_V)\\
  &H^q(\overline{Y}_{24},F_V)   \ar[r] \ar[ur]
             & H^q(\overline{Y}_{34},F_V) \ar[ur]
}
$$

The direct sum of the first column will be $E_1^{0,q}$; the direct
sum of the second column will be $E_1^{1,q}$; and $E_1^{2,q} =
H^q(\overline{Y}_B,F_V) $. We have non-zero terms when $q=2,3,4$
or $6$. Similarly, to the Mayer-Vietoris sequence, we want every
square at the $E_1$ level
 to be anti-commutative. It can be achieved in the following way.
First, consider the maps induced by the inclusion of the boundary
components. Then the squares will commute. Then change the sign of
every other arrow mapping a subspace of $E_1^{0,q}$ to a subspace
of $E_1^{1,q}$ as it is done in the definition of the spectral
sequence. Then the squares will anti-commute. \bth{6.1} The above
spectral sequence stabilizes at $E_2$ level. It converges to the
cohomology of the boundary of the Borel-Serre compactification
associated to $GL_4(\Z)$, which is
$$H^i_\d(GL_4(\Z),V)=\left\{
\begin{tabular}{lll}
$(0|0|0|n)\oplus (\overline{n-3,-1}|2|2)\oplus (\overline{n-3,1}|0|2)$
    & $i=3,$\\
$0$
    & $i\neq 3,$
\end{tabular}\right.
$$
where $$(a_1|a_2|\dots |a_k)=\otimes_{i=1}^k H^0(GL_1(\Z),L[a_i]),$$
and
$$(\overline{a_1,a_2}|a_3|a_4)=H^1_{cusp}(GL_2(\Z),L[a_1,a_2])
\otimes (a_3|a_4).$$
\eth
\proof
We consider all non-vanishing
terms of the spectral sequence at $E_1$ level. The non-vanishing
 terms occur
at $q=2,3,4$ and $6$. For a fixed $q$ we have arrows going in
 direction of the index $p$ induced by the inclusion of the
 parabolic subgroups. We compute
$kernel/image$ for these arrows in order to find the $E_2$ level
of the spectral sequence. As a consequence we find that the
spectral sequence degenerates at $E_2$ level. Then we compute the
cohomology to which it converges, which is the cohomology of the
boundary.
\subsection{Computation of $E_2^{*,2}$}
For the $E_1^{p,2}$-terms the only non-zero cohomologies come are
$H^2(P_{12},V)$ and $H^2(B,V)$. We have
$$(n-3,1|0|2) \rightarrow (0|n-2|0|2).$$
Therefore,
$$E_2^{p,2}=\left\{
\begin{tabular}{ll}
$(\overline{n-3,1}|0|2)$ & $p=1$\\
$0$                      & $p \neq1$
\end{tabular}
\right.
$$
\subsection{Computation of $E_2^{*,3}$}
First we consider the case $n>5$. Now we describe the $E_1^{*,3}$
terms. Consider the columns of the diagram below. Break each
column into pairs of vector spaces. Each pair comes one parabolic
subgroup. For example $(0|0|0|n)$ and $(\overline{n-3,-1}|2|2)$
come from  third cohomology of $P_{13}$. The two vector spaces
below come from the third cohomology of $P_{12,34}$. The maps
correspond to the inclusion of the parabolic subgroups.
$$
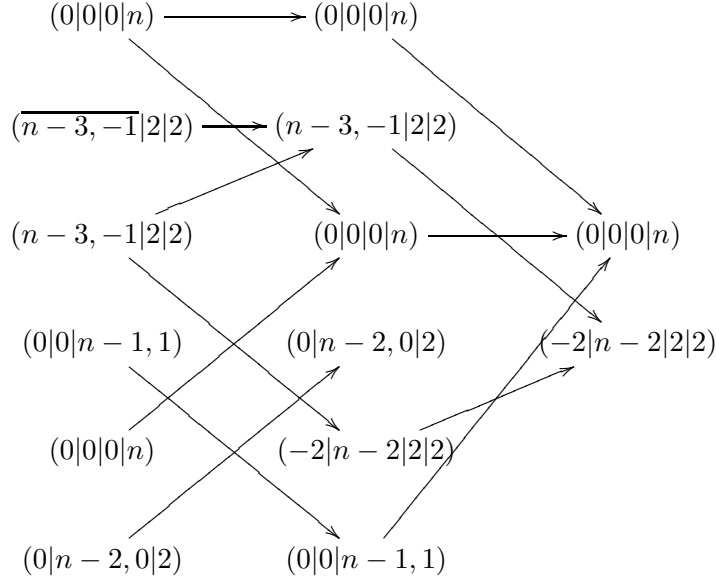
\xymatrix{
(0|0|0|n)               \ar[ddr]  \ar[r]  & (0|0|0|n)\ar[ddr]\\
(\overline{n-3,-1}|2|2) \ar[r]            & (n-3,-1|2|2)\ar[ddr]\\
(n-3,-1|2|2)            \ar[ur]  \ar[ddr] & (0|0|0|n)   \ar[r] & (0|0|0|n)\\
(0|0|n-1,1)             \ar[ddr]          & (0|n-2,0|2)        & (-2|n-2|2|2)\\
(0|0|0|n)               \ar[uur]          & (-2|n-2|2|2) \ar[ur]\\
(0|n-2,0|2)             \ar[uur]          & (0|0|n-1,1) \ar[uuur]
}
$$
There are many cancelation which occur when passing to $E_2$
level. In order to follow the cancelation one considers the
connected graph of the above diagram. There are 3 connected
graphs: one containing the space
 $(0|0|0|n)$ coming from the 3rd cohomology of the
Borel subgroup, and another containing $(-2|n-2|2|2)$ again from
the 3rd cohomology of the Borel subgroup, and the 3rd containing
$(0|n-2,0|2)$ from the 3rd cohomology of $P_{24}$. Consider the
graph containing $(0|0|0|n)$. The only term that is not cancelled
at $E_2$ level is the vector space $(0|0|0|n)$ which comes from
the parabolic group $P_{24}$. Now consider the second connected
graph, containing $(-2|n-2|2|2)$. After cancelation the only
vector space left is $(\overline{n-3,-1}|2|2)$ coming from
$P_{13}$. For the 3rd connected graph, there are two vertices
corresponding to  $(0|n-2,0|2)$. So they cancel and do not
contribute to the $E_2$ level. Thus, for $n>4$ we have
$$E_2^{p,3}=\left\{
\begin{tabular}{ll}
$(0|0|0|n)\oplus (\overline{n-3,-1}|2|2)$ & $p=0$,\\
$0$                                       & $p\neq 0.$
\end{tabular}
\right.
$$

Now we have to examine the case $n=4$. The vector spaces are all
the same as in the case $n>4$ except the exchange of
$(\overline{n-3,-1}|2|2)$ with $(0|n-2,0|2)$ in the 3rd cohomology
of $P_{13}$. Note also that for $n=4$, we have $(0|n-2,0|2)=0$.
Then the $E_1^{*,3}$ terms form the following anticommutative
diagram:
$$
\xymatrix{
(0|0|0|4)               \ar[ddr]  \ar[r]  & (0|0|0|4)\ar[ddr]\\
0                                         & (4-3,-1|2|2)\ar[ddr]\\
(4-3,-1|2|2)            \ar[ur]  \ar[ddr] & (0|0|0|4)   \ar[r] & (0|0|0|4)\\
(0|0|4-1,1)             \ar[ddr]          & 0                  & (-2|4-2|2|2)\\
(0|0|0|4)               \ar[uur]          & (-2|4-2|2|2) \ar[ur]\\
0                                         & (0|0|4-1,1) \ar[uuur]
}
$$


There are 2 connected graphs in the above diagram. One containing
the vector space $(0|0|0|4)$ coming from the Borel subgroup. The
other containing the vector space $(-2|4-2|2|2)$ again coming from
the Borel subgroup. Consider the graph containing $(0|0|0|4)$. The
only terms that is not canceled at $E_2$ level is the vector space
$(0|0|0|4)$ which comes from the parabolic group $P_{24}$. Now
consider the second connected graph, containing $(-2|4-2|2|2)$.
All of its terms of that graph cancel when passing to $E_2$ level.
Thus, for $w=4$ we have
$$E_2^{p,3}=\left\{
\begin{tabular}{ll}
$(0|0|0|4)$ & $p=0$,\\
$0$        & $p\neq 0.$
\end{tabular}
\right.
$$
\subsection{Computation of $E_2^{*,4}$}
For $q=4$ the only non-zero terms at $E_1$ level come from $P_{13}$ and
$P_{23}$. We have
$$H^4(P_{13},V))\rightarrow H^4(P_{23},V)).$$
From the first theorem (theorem 5.1) in the section "Cohomology
of the parabolic subgroups of $GL_4$" we obtain
$$(-2|n-2,2|2)\rightarrow (-2|n-2,2|2).$$
Therefore,
$$E_2^{*,4}=0.$$
\subsection{Computation of $E_2^{*,6}$}
When $q=6$, for all even $w$, the non-zero terms give
$$E_1^{*,6}: H^6(P_{24},V)\rightarrow  H^6(P_{34},V)\rightarrow  H^6(B,V),$$
which are isomorphic to
$$(-2|0|\overline{w-1,3}) \rightarrow  (-2|0|w-1,3) \rightarrow
(-2|0|2|w)$$
from theorem 5.1.
The above sequence is exact. Therefore,
$$E_2^{*,6} =0.$$

The spectral sequence degenerates at $E_2$ level. Therefore, we can find what
is the cohomology of the boundary of the Borel-Serre compactification
associated to $GL_4(\Z)$ with coefficients in the sheaf $F_V$ associated to
$$V=S^{n-4}V_4\otimes det.$$

Let us recall the notation that we are going to use.
By $H^i_\d(GL_4(\Z),V)$ we mean the cohomology of the boundary of
the Borel-Serre compactification associated to $GL_4(\Z)$ with coefficients
the sheaf $F_V$
For even $n$ greater than $4$ we have
$$H^i_\d(GL_4(\Z),S^{n-4}V_4\otimes det)=\left\{
\begin{tabular}{ll}
$(0|0|0|n)\oplus (\overline{n-3,-1}|2|2) \oplus (\overline{n-3,1}|0|2)$
    & $i=3$,\\
$0$ & $i\neq 3$.
\end{tabular}
\right.
$$
Note that the first two summands for the 3rd cohomology of the boundary
come from 3rd cohomology of the maximal parabolic subgroups. And the last
summand comes from the 2nd cohomology of a non-maximal parabolic subgroup.
Since it comes from second cohomology of a parabolic subgroup,
but it contributes in the 3rd cohomology of the boundary, it is called a ghost
class.
\sectionnew{Cohomology of $GL_4(\Z)$}
We are going to show that the ghost class do not enter in the
Eisenstein cohomology of $GL_4(\Z)$ which coinsides with the whole
cohomology of $GL_4(\Z)$.

Since the cohomology of the boundary is concentrated in degree 3,
it is enough to compute homological Euler characteristic of $GL_4(\Z)$
 with coefficients in $S^{n-4}V_4\otimes det$. Recall the homological
Euler characteristic of an arithmetic group $\Gamma$ with
coefficients in a finite dimensional representation is
$$\chi_h(\Gamma,V)=\sum_i (-1)^i\dim H^i(\Gamma,V).$$
Note that $S^{n-4}V_4\otimes det=L[n-3,1,1,1]$ and
 $S^{n-2}V_2\otimes det=L[n-1,1]=L[n-3,-1]$
Form \cite{???} we know that
$$\chi_h(GL_4(\Z),S^{n-4}V_4\otimes det)=
\chi_h(GL_2(\Z),S^{n-2}V_2\otimes det).$$
Therefore, for even $n$ greater than $4$, we have
$$H^i(GL_4(\Z),S^{n-4}V_4\otimes det)=\left\{
\begin{tabular}{ll}
$(0|0|0|n)\oplus (\overline{n-3,-1}|2|2)$
    & $i=3$,\\
$0$ & $i\neq 3$.
\end{tabular}
\right.
$$
In the case $n=4$ we use the same argument.
$$H^i_\d(GL_4(\Z),det)=\left\{
\begin{tabular}{ll}
$(0|0|0|4)$  & $i=3$,\\
$0$          & $i\neq 3$.
\end{tabular}
\right.
$$
Also, the homological Euler characteristic gives
$$\chi_h(GL_4(\Z),det)=-1.$$
Therefore, for $n=4$ the cohomology of the boundary coincides with
the Eisenstein cohomology. And we have
$$H^i_{Eis}(GL_4(\Z),det)=\left\{
\begin{tabular}{ll}
$(0|0|0|4)$  & $i=3$,\\
$0$          & $i\neq 3$.
\end{tabular}
\right.
$$
On the other hand,
$$H^i_{cusp}(SL_4(\Z),\Q)=0.$$
Therefore,
$$H^i_{cusp}(GL_4(\Z),det)=0.$$
And we conclude that
$$H^i(GL_4(\Z),det)=\left\{
\begin{tabular}{ll}
$(0|0|0|4)$  & $i=3$,\\
$0$          & $i\neq 3$.
\end{tabular}
\right.
$$

\renewcommand{\em}{\textrm}

\begin{small}

\renewcommand{\refname}{ {\flushleft\normalsize\bf{References}} }
    
\end{small}

\end{document}